\documentclass{amsart}
\usepackage[dvips]{graphicx}

\newtheorem{thm}{Theorem}

\theoremstyle{definition}
\newtheorem{df}[thm]{\sc Definition}
\newtheorem{remark}[thm]{\sc Remark}
\newtheorem{example}[thm]{\sc Example}
\newtheorem*{acknowledgment}{\sc Acknowledgment}

\numberwithin{equation}{section}

\begin{document}

\title[Ehrhart polynomials of 3-dimensional simple polytopes]
{Ehrhart polynomials of 3-dimensional simple integral convex polytopes}

\author{Yusuke Suyama}
\address{Department of Mathematics, Graduate School of Science, Osaka City University,
3-3-138 Sugimoto, Sumiyoshi-ku, Osaka 558-8585 JAPAN}
\email{d15san0w03@st.osaka-cu.ac.jp}
\thanks{}



\date{\today}

\dedicatory{}

\begin{abstract}
We give an explicit formula on the Ehrhart polynomial
of a 3-dimensional simple integral convex polytope by using toric geometry.
\end{abstract}

\maketitle

\section{Introduction}

Let $P \subset \mathbb{R}^d$ be an integral convex polytope of dimension $d$,
that is, a convex polytope whose vertices have integer coordinates.
For a non-negative integer $l$, we write $lP=\{lx \mid x \in P\}$.
Ehrhart \cite{Ehrhart} proved that the number of lattice points in $lP$
can be expressed by a polynomial in $l$ of degree $d$:
\begin{equation*}
|(lP) \cap \mathbb{Z}^d|=c_dl^d+c_{d-1}l^{d-1}+\cdots+c_0.
\end{equation*}
This polynomial is called the {\it Ehrhart polynomial} of $P$.
It is known that:
\begin{enumerate}
\item $c_0=1$.
\item $c_{d-1}$ is half of the sum of relative volumes of facets of $P$
(\cite[Theorem 5.6]{BR}).
\item $c_d$ is the volume of $P$ (\cite[Corollary 3.20]{BR}).
\end{enumerate}
However, we have no formula on other coefficients of Ehrhart polynomials.
In particular, we do not know a formula on $c_1$
for a general 3-dimensional integral convex polytope.
In this paper, we find an explicit formula on $c_1$ of the Ehrhart polynomial
of a 3-dimensional {\it simple} integral convex polytope, see Theorem \ref{main}.

Pommersheim \cite{Pommersheim} gave a method for computing the $(d-2)$-nd coefficient
of the Ehrhart polynomial of a $d$-dimensional simple integral convex polytope $P$
by using toric geometry.
He obtained an explicit description of the Ehrhart polynomial of a tetrahedron
by using this method.
Our formula is obtained by using this method
for a general 3-dimensional simple integral convex polytope.

The structure of the paper is as follows.
In Section 2, we state the main theorem and give a few examples.
In Section 3, we give a proof of the main theorem.

\begin{acknowledgment}
This work was supported by Grant-in-Aid for JSPS Fellows 15J01000.
The author wishes to thank his supervisor, Professor Mikiya Masuda,
for his continuing support.
\end{acknowledgment}

\section{The main theorem}

Let $P \subset \mathbb{R}^3$ be a 3-dimensional simple integral convex polytope,
and let $F_1, \ldots, F_n$ be the facets of $P$.
For $k=1, \ldots, n$, we denote by $v_k \in \mathbb{Z}^3$
the inward-pointing primitive normal vector of $F_k$.
For an edge $E$ of $P$, we denote by $\mathrm{Vol}(E)$
the relative volume of $E$, that is, the length of $E$ measured with respect to
the lattice of rank one in the line containing $E$.

\begin{df}
For each edge $E=F_{k_1} \cap F_{k_2}$ of $P$,
we define an integer $m(E)$ and a rational number $s(E)$ as follows:
\begin{enumerate}\setlength{\itemsep}{-1mm}
\item We define $m(E)=|((\mathbb{R}v_{k_1}+\mathbb{R}v_{k_2}) \cap \mathbb{Z}^3)/
(\mathbb{Z}v_{k_1}+\mathbb{Z}v_{k_2})|$. \\
\item There exists a basis $e_1, e_2$ for
$(\mathbb{R}v_{k_1}+\mathbb{R}v_{k_2}) \cap \mathbb{Z}^3$
such that $v_{k_1}=e_1$ and $v_{k_2}=pe_1+qe_2$ for some $q>p \geq 0$.
Then we define $s(E)=s(p, q)$, where $s(p, q)$ is the Dedekind sum, which is defined by
\begin{equation*}
s(p, q)
=\sum_{i=1}^q \left(\left(\frac{i}{q}\right)\right)\left(\left(\frac{pi}{q}\right)\right),\quad
((x))=\left\{\begin{array}{ll}
x-[x]-\frac{1}{2} & (x \notin \mathbb{Z}),\\
0 & (x \in \mathbb{Z}).
\end{array}\right.
\end{equation*}
\end{enumerate}
\end{df}

\begin{remark}
We have $q=m(E)$. Although $p$ is not uniquely determined,
$s(p, q)$ does not depend on the choice of $e_1, e_2$.
Thus $s(E)$ is well-defined.
\end{remark}

\begin{df}\label{C}
For each facet $F$ of $P$, we define a rational number $C(F)$ as follows.
We name vertices and facets around $F$ as in Figure \ref{Ehrhart}.
We denote by $v \in \mathbb{Z}^3$
the inward-pointing primitive normal vector of $F$.

\begin{figure}[htbp]
\begin{center}
\includegraphics[width=6.5cm]{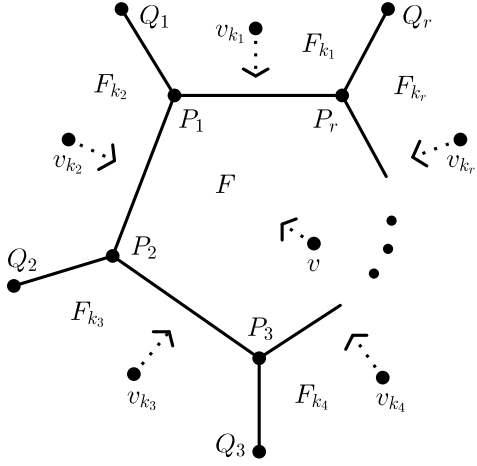}
\caption{vertices and facets around $F$.}
\label{Ehrhart}
\end{center}
\end{figure}

For $i=1, \ldots, r$, we define
\begin{equation*}
\varepsilon_i=\mathrm{det}(v, v_{k_{i+1}}, v_{k_i})>0,\quad
a_i=\cfrac{\langle \overrightarrow{P_{i-1}Q_{i-1}}, v_{k_{i+1}} \rangle}
{\varepsilon_i \langle \overrightarrow{P_{i-1}Q_{i-1}}, v \rangle},\quad
b_i=\cfrac{\langle \overrightarrow{P_iP_{i+1}}, v_{k_{i-1}} \rangle}
{\varepsilon_{i-1} \langle \overrightarrow{P_iP_{i+1}}, v_{k_i} \rangle},
\end{equation*}
where
$v_{k_0}=v_{k_r}, v_{k_{r+1}}=v_{k_1}, \varepsilon_0=\varepsilon_r, P_0=P_r, P_{r+1}=P_1, Q_0=Q_r$
and $\langle \cdot, \cdot \rangle$ is the standard inner product on $\mathbb{R}^3$.
Then we define
\begin{equation*}
C(F)=-\sum_{2 \leq i<j \leq r}a_i
\left|\begin{array}{ccccc}
b_{i+1} & \varepsilon_{i+1}^{-1} & 0 & \cdots & 0 \\
\varepsilon_{i+1}^{-1} & b_{i+2} & \varepsilon_{i+2}^{-1} & \ddots & \vdots \\
0 & \varepsilon_{i+2}^{-1} & \ddots & \ddots & 0 \\
\vdots & \ddots & \ddots & b_{j-2} & \varepsilon_{j-2}^{-1} \\
0 & \cdots & 0 & \varepsilon_{j-2}^{-1} & b_{j-1} \\
\end{array}\right|
\varepsilon_i \varepsilon_{i+1} \cdots \varepsilon_{j-1}
\frac{\mathrm{Vol}(P_{j-1}P_j)}{m(P_{j-1}P_j)},
\end{equation*}
where $P_{j-1}P_j$ is the edge whose endpoints are $P_{j-1}$ and $P_j$,
and the determinants above are understood to be one when $j=i+1$.
\end{df}

\begin{remark}
The proof of Theorem \ref{main} below shows that
$C(F)$ does not depend on the choice of $F_{k_1}$.
\end{remark}

The following is our main theorem:

\begin{thm}\label{main}
Let $P \subset \mathbb{R}^3$ be a 3-dimensional simple integral convex polytope,
and let $E_1, \ldots, E_m$ and $F_1, \ldots, F_n$
be the edges and the facets of $P$, respectively.
Then the coefficient $c_1$ of the Ehrhart polynomial
$|(lP) \cap \mathbb{Z}^3|=c_3l^3+c_2l^2+c_1l+c_0$ is given by
\begin{equation*}
\sum_{j=1}^m \left(s(E_j)+\frac{1}{4}\right)\mathrm{Vol}(E_j)
+\frac{1}{12}\sum_{k=1}^n C(F_k).
\end{equation*}
\end{thm}

\begin{example}
Let $a, b, c$ be positive integers with $\mathrm{gcd}(a, b, c)=1$
and let $P \subset \mathbb{R}^3$ be the tetrahedron with vertices
\begin{equation*}
O=\left(\begin{array}{c}0\\0\\0\end{array}\right),\quad
P_1=\left(\begin{array}{c}a\\0\\0\end{array}\right),\quad
P_2=\left(\begin{array}{c}0\\b\\0\end{array}\right),\quad
P_3=\left(\begin{array}{c}0\\0\\c\end{array}\right).
\end{equation*}
We put $A=\mathrm{gcd}(b,c), B=\mathrm{gcd}(a,c), C=\mathrm{gcd}(a,b)$
and $d=ABC$. Then we have the following table:
\begin{table}[htbp]
\begin{center}
\begin{tabular}{|c||c|c|c|c|c|c|}
\hline
edge $E$ & $OP_1$ & $OP_2$ & $OP_3$ & $P_1P_2$ & $P_1P_3$ & $P_2P_3$ \\
\hline
$\mathrm{Vol}(E)$ & $a$ & $b$ & $c$ & $C$ & $B$ & $A$ \\
\hline
$m(E)$ & $1$ & $1$ & $1$ & $cC/d$ & $bB/d$ & $aA/d$ \\
\hline
$s(E)$ & $0$ & $0$ & $0$ &
$-s\left(\cfrac{ab}{d}, \cfrac{cC}{d}\right)$ &
$-s\left(\cfrac{ac}{d}, \cfrac{bB}{d}\right)$ &
$-s\left(\cfrac{bc}{d}, \cfrac{aA}{d}\right)$ \\
\hline
\end{tabular}
\begin{tabular}{|c||c|c|c|c|}
\hline
facet $F$ & $OP_1P_2$ & $OP_1P_3$ & $OP_2P_3$ & $P_1P_2P_3$ \\
\hline
\shortstack{inward-pointing primitive\\normal vector of $F$} &
$\left(\begin{array}{c}0\\0\\1\end{array}\right)$ &
$\left(\begin{array}{c}0\\1\\0\end{array}\right)$ &
$\left(\begin{array}{c}1\\0\\0\end{array}\right)$ &
$\left(\begin{array}{c}-bc/d\\-ac/d\\-ab/d\end{array}\right)$ \\
\hline
$C(F)$ & $ab/c$ & $ac/b$ &$bc/a$ & $d^2/(abc)$ \\
\hline
\end{tabular}
\caption{the values of $\mathrm{Vol}(E), s(E)$ and $C(F)$.}
\label{values1}
\end{center}
\end{table}

Thus we have
\begin{align*}
&\sum_{E:\mathrm{edge}} \left(s(E)+\frac{1}{4}\right)\mathrm{Vol}(E)
+\frac{1}{12}\sum_{F:\mathrm{facet}} C(F)\\
&=\frac{a}{4}+\frac{b}{4}+\frac{c}{4}
+\left(-s\left(\frac{ab}{d}, \frac{cC}{d}\right)+\frac{1}{4}\right)C
+\left(-s\left(\frac{ac}{d}, \frac{bB}{d}\right)+\frac{1}{4}\right)B\\
&+\left(-s\left(\frac{bc}{d}, \frac{aA}{d}\right)+\frac{1}{4}\right)A
+\frac{1}{12}\left(\frac{ab}{c}+\frac{ac}{b}+\frac{bc}{a}+\frac{d^2}{abc}\right),
\end{align*}
which coincides with the formula in \cite[Theorem 5]{Pommersheim}.
\end{example}

\begin{example}
Let $a$ and $c$ be positive integers and $b$ be a non-negative integer.
Consider the convex hull $P \subset \mathbb{R}^3$ of the six points
\begin{align*}
&O=\left(\begin{array}{c}0\\0\\0\end{array}\right),\quad
A=\left(\begin{array}{c}a\\0\\0\end{array}\right),\quad
B=\left(\begin{array}{c}0\\a\\0\end{array}\right),\\
&O'=\left(\begin{array}{c}b\\0\\c\end{array}\right),\quad
A'=\left(\begin{array}{c}a+b\\0\\c\end{array}\right),\quad
B'=\left(\begin{array}{c}b\\a\\c\end{array}\right).
\end{align*}
$P$ is a 3-dimensional simple polytope.
We put $g=\mathrm{gcd}(b,c)$. Then we have the following table:
\begin{table}[htbp]
\begin{center}
\footnotesize\begin{tabular}{|c||c|c|c|c|c|c|c|c|c|}
\hline
edge $E$ & $OA$ & $OB$ & $AB$ & $OO'$ & $AA'$ & $BB'$ & $O'A'$ & $O'B'$ & $A'B'$ \\
\hline
$\mathrm{Vol}(E)$ & $a$ & $a$ & $a$ & $g$ & $g$ & $g$ & $a$ & $a$ & $a$ \\
\hline
$m(E)$ & $1$ & $c/g$ & $c/g$ & $1$ & $1$ & $c/g$ & $1$ & $c/g$ & $c/g$ \\
\hline
$s(E)$ & $0$ & $-s\left(\frac{b}{g}, \frac{c}{g}\right)$ & $s\left(\frac{b}{g}, \frac{c}{g}\right)$ &
$0$ & $0$ & $-s\left(1, \frac{c}{g}\right)$ &
$0$ & $s\left(\frac{b}{g}, \frac{c}{g}\right)$ & $-s\left(\frac{b}{g}, \frac{c}{g}\right)$ \\
\hline
\end{tabular}
\begin{tabular}{|c||c|c|c|c|c|}
\hline
facet $F$ & $OAB$ & $OAA'O'$ & $OBB'O'$ & $ABB'A'$ & $O'A'B'$ \\
\hline
\shortstack{inward-pointing primitive\\normal vector of $F$} &
$\left(\begin{array}{c}0\\0\\1\end{array}\right)$ &
$\left(\begin{array}{c}0\\1\\0\end{array}\right)$ &
$\left(\begin{array}{c}c/g\\0\\-b/g\end{array}\right)$ &
$\left(\begin{array}{c}-c/g\\-c/g\\b/g\end{array}\right)$ &
$\left(\begin{array}{c}0\\0\\-1\end{array}\right)$ \\
\hline
$C(F)$ & $0$ & $c$ & $g^2/c$ & $g^2/c$ & $0$ \\
\hline
\end{tabular}
\caption{the values of $\mathrm{Vol}(E), s(E)$ and $C(F)$.}
\label{values2}
\end{center}
\end{table}
\begin{figure}[htbp]
\begin{center}
\includegraphics[width=7cm]{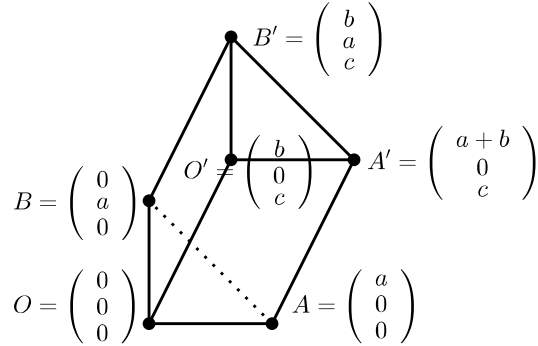}
\caption{the simple polytope $P$.}
\label{example2}
\end{center}
\end{figure}

Thus we have
\begin{align*}
&\sum_{E:\mathrm{edge}} \left(s(E)+\frac{1}{4}\right)\mathrm{Vol}(E)
+\frac{1}{12}\sum_{F:\mathrm{facet}} C(F)\\
&=-s\left(1, \frac{c}{g}\right)g+\frac{3a}{2}+\frac{3g}{4}
+\frac{1}{12}\left(c+\frac{2g^2}{c}\right)\\
&=-g\sum_{i=1}^{c/g-1}\left(\frac{i}{\frac{c}{g}}-\frac{1}{2}\right)^2
+\frac{3a}{2}+\frac{3g}{4}+\frac{c}{12}+\frac{g^2}{6c}\\
&=-g\sum_{i=1}^{c/g-1}\left(\frac{g^2}{c^2}i^2-\frac{g}{c}i+\frac{1}{4}\right)
+\frac{3a}{2}+\frac{3g}{4}+\frac{c}{12}+\frac{g^2}{6c}\\
&=-\frac{g^3}{c^2}\frac{\left(\frac{c}{g}-1\right)\frac{c}{g}\left(\frac{2c}{g}-1\right)}{6}
+\frac{g^2}{c}\frac{\left(\frac{c}{g}-1\right)\frac{c}{g}}{2}-g\frac{\frac{c}{g}-1}{4}
+\frac{3a}{2}+\frac{3g}{4}+\frac{c}{12}+\frac{g^2}{6c}\\
&=\frac{3a}{2}+g.
\end{align*}
On the other hand, since
\begin{equation*}
\#\{(x, y) \in \mathbb{Z}^2 \mid (x, y, z) \in lP\}
=\left\{\begin{array}{ll}
\cfrac{(al+1)(al+2)}{2} & ((c/g)|z), \\
\cfrac{al(al+1)}{2} & ((c/g)\mid \hspace{-.67em}/z) \end{array}\right.
\end{equation*}
for $z=0, 1, \ldots, cl$, we have
\begin{align*}
|(lP) \cap \mathbb{Z}^3|&=\frac{(al+1)(al+2)}{2}(gl+1)+\cfrac{al(al+1)}{2}((cl+1)-(gl+1))\\
&=\frac{a^2c}{2}l^3+\frac{1}{2}\left(a^2+ac+2ag\right)l^2+\left(\frac{3a}{2}+g\right)l+1.
\end{align*}
The coefficient of $l$ is also $3a/2+g$.
\end{example}

\section{Proof of Theorem \ref{main}}

First we recall some facts about toric geometry, see \cite{Fulton} for details.
Let $P \subset \mathbb{R}^d$ be a $d$-dimensional integral convex polytope.
We define a cone
\begin{equation*}
\sigma_F=\{v \in \mathbb{R}^d \mid
\langle u'-u, v\rangle \geq 0\ \forall u' \in P, \forall u \in F\}
\end{equation*}
for each face $F$ of $P$.
Then the set
\begin{equation*}
\Delta_P=\{\sigma_F \mid F\mbox{ is a face of }P\}
\end{equation*}
of such cones
forms a fan in $\mathbb{R}^d$, which is called the {\it normal fan} of $P$.
Let $X(\Delta_P)$ be the associated projective toric variety.
We denote by $V(\sigma)$ the subvariety of $X(\Delta_P)$
corresponding to $\sigma \in \Delta_P$.
Let $\mathrm{Td}_i(X(\Delta_P)) \in A_i(X(\Delta_P))_\mathbb{Q}$
be the $i$-th Todd class in the Chow group of $i$-cycles with rational coefficients.

\begin{thm}\label{Todd}
Let $P \subset \mathbb{R}^d$ be a $d$-dimensional integral convex polytope
and $|(lP) \cap \mathbb{Z}^d|=c_dl^d+c_{d-1}l^{d-1}+\cdots+c_0$
be its Ehrhart polynomial.
If $\mathrm{Td}_i(X(\Delta_P))$ has an expression of the form
$\sum_Fr_F[V(\sigma_F)]$ with $r_F \in \mathbb{Q}$,
then we have $c_i=\sum_Fr_F\mathrm{Vol}(F)$,
where $[V(\sigma_F)]$ is the class of $V(\sigma_F)$ in the Chow group
and $\mathrm{Vol}(F)$ is the relative volume of $F$.
\end{thm}

Now we assume that $d=3$ and $P$ is simple.
Then the associated toric variety $X(\Delta_P)$ is $\mathbb{Q}$-factorial
and we know the ring structure of the Chow ring $A^*(X(\Delta_P))_\mathbb{Q}$
with rational coefficients.
Let $E_1, \ldots, E_m$ and $F_1, \ldots, F_n$
be the edges and the facets of $P$, respectively. We have
\begin{equation}\label{1}
\sum_{k=1}^n\langle u, v_k\rangle[V(\sigma_{F_k})]=0 \quad \forall u \in (\mathbb{Q}^3)^*.
\end{equation}
If $F_{k_1}$ and $F_{k_2}$ are distinct, then
\begin{equation}\label{2}
[V(\sigma_{F_{k_1}})][V(\sigma_{F_{k_2}})]
=\left\{\begin{array}{ll}
\frac{1}{m(E_j)}[V(\sigma_{E_j})] & (1 \leq \exists j \leq m: F_{k_1} \cap F_{k_2}=E_j), \\
0 & (F_{k_1} \cap F_{k_2}=\emptyset)
\end{array}\right.
\end{equation}
in $A^*(X(\Delta_P))_\mathbb{Q}$.

Pommersheim gave an expression of $\mathrm{Td}_{d-2}(X(\Delta_P))$
for a $d$-dimensional simple integral convex polytope $P \subset \mathbb{R}^d$.
In the case where $d=3$, we have the following:

\begin{thm}[Pommersheim \cite{Pommersheim}]\label{P}
If $P \subset \mathbb{R}^3$ is a 3-dimensional simple integral convex polytope, then
\begin{equation*}
\mathrm{Td}_1(X(\Delta_P))=\sum_{j=1}^m \left(s(E_j)+\frac{1}{4}\right)[V(\sigma_{E_j})]
+\frac{1}{12}\sum_{k=1}^n [V(\sigma_{F_k})]^2.
\end{equation*}
\end{thm}

We use the notation in Definition \ref{C}.
It suffices to show
\begin{equation*}
[V(\sigma_{F})]^2=-\sum_{2 \leq i<j \leq r}a_i
\left|\begin{array}{ccccc}
b_{i+1} & \varepsilon_{i+1}^{-1} & 0 & \cdots & 0 \\
\varepsilon_{i+1}^{-1} & b_{i+2} & \varepsilon_{i+2}^{-1} & \ddots & \vdots \\
0 & \varepsilon_{i+2}^{-1} & \ddots & \ddots & 0 \\
\vdots & \ddots & \ddots & b_{j-2} & \varepsilon_{j-2}^{-1} \\
0 & \cdots & 0 & \varepsilon_{j-2}^{-1} & b_{j-1} \\
\end{array}\right|
\frac{\varepsilon_i \varepsilon_{i+1} \cdots \varepsilon_{j-1}}{m(P_{j-1}P_j)}
[V(\sigma_{P_{j-1}P_j})]
\end{equation*}
for each facet $F$ of $P$.

We put
\begin{equation*}
D(s, t)=\left|\begin{array}{ccccc}
b_s & \varepsilon_s^{-1} & 0 & \cdots & 0 \\
\varepsilon_s^{-1} & b_{s+1} & \varepsilon_{s+1}^{-1} & \ddots & \vdots \\
0 & \varepsilon_{s+1}^{-1} & \ddots & \ddots & 0 \\
\vdots & \ddots & \ddots & b_{t-1} & \varepsilon_{t-1}^{-1} \\
0 & \cdots & 0 & \varepsilon_{t-1}^{-1} & b_t \\
\end{array}\right|
\end{equation*}
for $2<s \leq t<r$ and $D(s, t)=1$ for $s>t$.
Define $u \in (\mathbb{Q}^3)^*$ by
$\langle u, v \rangle=1, \langle u, v_{k_1} \rangle=0, \langle u, v_{k_2} \rangle=0$.
By (\ref{1}) and (\ref{2}), we have
\begin{equation*}
[V(\sigma_{F})]^2=-[V(\sigma_{F})]\sum_{j=1}^r\langle u, v_{k_j}\rangle[V(\sigma_{F_{k_j}})]
=-\sum_{j=3}^r\frac{\langle u, v_{k_j}\rangle}{m(P_{j-1}P_j)}[V(\sigma_{P_{j-1}P_j})].
\end{equation*}
Hence it suffices to show
\begin{equation}\label{induction}
\langle u, v_{k_j}\rangle=\sum_{i=2}^{j-1}a_i
D(i+1, j-1)\varepsilon_i \varepsilon_{i+1} \cdots \varepsilon_{j-1}
\end{equation}
for any $j=3, \ldots, r$.

First we claim that
\begin{equation}\label{relation}
\varepsilon_{j-1}^{-1}v_{k_{j-1}}+\varepsilon_j^{-1}v_{k_{j+1}}=a_jv+b_jv_{k_j}
\end{equation}
for any $j=2, \ldots, r-1$. By Cramer's rule, we have
\begin{align*}
v_{k_{j+1}}
&=\frac{\mathrm{det}(v_{k_{j+1}}, v_{k_j}, v_{k_{j-1}})}{\mathrm{det}(v, v_{k_j}, v_{k_{j-1}})}v
+\frac{\mathrm{det}(v, v_{k_{j+1}}, v_{k_{j-1}})}{\mathrm{det}(v, v_{k_j}, v_{k_{j-1}})}v_{k_j}
+\frac{\mathrm{det}(v, v_{k_j}, v_{k_{j+1}})}{\mathrm{det}(v, v_{k_j}, v_{k_{j-1}})}v_{k_{j-1}}\\
&=\frac{\mathrm{det}(v_{k_{j+1}}, v_{k_j}, v_{k_{j-1}})}{\varepsilon_{j-1}}v
+\frac{\mathrm{det}(v, v_{k_{j+1}}, v_{k_{j-1}})}{\varepsilon_{j-1}}v_{k_j}
-\frac{\varepsilon_j}{\varepsilon_{j-1}}v_{k_{j-1}}.
\end{align*}
So we have
\begin{align}\label{relation2}
&\begin{aligned}
&\varepsilon_{j-1}^{-1}v_{k_{j-1}}+\varepsilon_j^{-1}v_{k_{j+1}}\\
&=\varepsilon_{j-1}^{-1}\varepsilon_j^{-1}\mathrm{det}(v_{k_{j+1}}, v_{k_j}, v_{k_{j-1}})v
+\varepsilon_{j-1}^{-1}\varepsilon_j^{-1}\mathrm{det}(v, v_{k_{j+1}}, v_{k_{j-1}})v_{k_j}.
\end{aligned}
\end{align}
Taking the inner product of both sides of (\ref{relation2})
with $\overrightarrow{P_{j-1}Q_{j-1}}$ gives
\begin{equation*}
\varepsilon_j^{-1}\langle \overrightarrow{P_{j-1}Q_{j-1}}, v_{k_{j+1}} \rangle
=\varepsilon_{j-1}^{-1}\varepsilon_j^{-1}\mathrm{det}(v_{k_{j+1}}, v_{k_j}, v_{k_{j-1}})
\langle \overrightarrow{P_{j-1}Q_{j-1}}, v \rangle,
\end{equation*}
which means
$a_j=\varepsilon_{j-1}^{-1}\varepsilon_j^{-1}\mathrm{det}(v_{k_{j+1}}, v_{k_j}, v_{k_{j-1}})$.
Taking the inner product of both sides of (\ref{relation2})
with $\overrightarrow{P_jP_{j+1}}$ gives
\begin{equation*}
\varepsilon_{j-1}^{-1}\langle \overrightarrow{P_jP_{j+1}}, v_{k_{j-1}} \rangle
=\varepsilon_{j-1}^{-1}\varepsilon_j^{-1}\mathrm{det}(v, v_{k_{j+1}}, v_{k_{j-1}})
\langle \overrightarrow{P_jP_{j+1}}, v_{k_j} \rangle,
\end{equation*}
which means
$b_j=\varepsilon_{j-1}^{-1}\varepsilon_j^{-1}\mathrm{det}(v, v_{k_{j+1}}, v_{k_{j-1}})$.
Thus (\ref{relation}) follows.

We show (\ref{induction}) by induction on $j$.
If $j=3$, then both sides are $a_2\varepsilon_2$.
If $j=4$, then both sides are $a_2b_3\varepsilon_2\varepsilon_3+a_3\varepsilon_3$.
Suppose $4 \leq j \leq r-1$. By (\ref{relation}) and the hypothesis of induction, we have
\begin{align*}
\langle u, v_{k_{j+1}}\rangle
&=\langle u, a_j\varepsilon_jv+b_j\varepsilon_jv_{k_j}
-\varepsilon_{j-1}^{-1}\varepsilon_jv_{k_{j-1}}\rangle\\
&=a_j\varepsilon_j+b_j\varepsilon_j\langle u, v_{k_j}\rangle
-\varepsilon_{j-1}^{-1}\varepsilon_j\langle u, v_{k_{j-1}}\rangle\\
&=a_j\varepsilon_j+b_j\varepsilon_j\sum_{i=2}^{j-1}a_i
D(i+1, j-1)\varepsilon_i \varepsilon_{i+1} \cdots \varepsilon_{j-1}\\
&-\varepsilon_{j-1}^{-1}\varepsilon_j\sum_{i=2}^{j-2}a_i
D(i+1, j-2)\varepsilon_i \varepsilon_{i+1} \cdots \varepsilon_{j-2}.
\end{align*}
On the other hand,
\begin{align*}
&\sum_{i=2}^ja_i
D(i+1, j)\varepsilon_i \varepsilon_{i+1} \cdots \varepsilon_j\\
&=a_j\varepsilon_j+a_{j-1}b_j\varepsilon_{j-1}\varepsilon_j
+\sum_{i=2}^{j-2}a_i
D(i+1, j)\varepsilon_i \varepsilon_{i+1} \cdots \varepsilon_j.
\end{align*}
Since
\begin{align*}
&\sum_{i=2}^{j-2}a_i
D(i+1, j)\varepsilon_i \varepsilon_{i+1} \cdots \varepsilon_j\\
&=\sum_{i=2}^{j-2}a_i(b_jD(i+1, j-1)-\varepsilon_{j-1}^{-2}D(i+1, j-2))
\varepsilon_i \varepsilon_{i+1} \cdots \varepsilon_j\\
&=b_j\varepsilon_j\sum_{i=2}^{j-2}
a_iD(i+1, j-1)\varepsilon_i \varepsilon_{i+1} \cdots \varepsilon_{j-1}\\
&-\varepsilon_{j-1}^{-1}\varepsilon_j\sum_{i=2}^{j-2}a_i
D(i+1, j-2)\varepsilon_i \varepsilon_{i+1} \cdots \varepsilon_{j-2},
\end{align*}
we have
\begin{align*}
&\sum_{i=2}^ja_i
D(i+1, j)\varepsilon_i \varepsilon_{i+1} \cdots \varepsilon_j\\
&=a_j\varepsilon_j+a_{j-1}b_j\varepsilon_{j-1}\varepsilon_j
+b_j\varepsilon_j\sum_{i=2}^{j-2}
a_iD(i+1, j-1)\varepsilon_i \varepsilon_{i+1} \cdots \varepsilon_{j-1}\\
&-\varepsilon_{j-1}^{-1}\varepsilon_j\sum_{i=2}^{j-2}a_i
D(i+1, j-2)\varepsilon_i \varepsilon_{i+1} \cdots \varepsilon_{j-2}\\
&=a_j\varepsilon_j+b_j\varepsilon_j\sum_{i=2}^{j-1}a_i
D(i+1, j-1)\varepsilon_i \varepsilon_{i+1} \cdots \varepsilon_{j-1}\\
&-\varepsilon_{j-1}^{-1}\varepsilon_j\sum_{i=2}^{j-2}a_i
D(i+1, j-2)\varepsilon_i \varepsilon_{i+1} \cdots \varepsilon_{j-2}\\
&=\langle u, v_{k_{j+1}}\rangle.
\end{align*}
Thus (\ref{induction}) holds for $j+1$. This completes the proof of Theorem \ref{main}.

\end{document}